\documentclass[12pt]{amsart}
\usepackage{amsmath,amssymb,amscd,amsfonts}
\newtheorem{prop}{Proposition}[section]
\newtheorem{thm}[prop]{Theorem}
\newtheorem{cor}[prop]{Corollary}

\newtheorem{lem}[prop]{Lemma}
\def\demo#1{\medskip\noindent{\bf#1\enspace}}

\def \calS {{\mathcal S}}
\def \calT {{\mathcal T}}

        \def\red{{\mathbb R}}
        
	\def\qed{{\mathbb Q}}
        
\begin{document}

\title{When size matters: subshifts and their related tiling spaces}

\author{Alex Clark and Lorenzo Sadun}

\address{Alex Clark: Department of Mathematics, University of North Texas,
             Denton, Texas 76203}
\email{AlexC@unt.edu}

\address{Lorenzo Sadun: Department of Mathematics, The University of Texas at Austin,
 Austin, TX 78712-1082 U.S.A.}
\email{sadun@math.utexas.edu}   



\begin{abstract}
We investigate the dynamics of substitution subshifts and their 
associated tiling spaces.  For a given subshift, the associated tiling
spaces are all homeomorphic, but their dynamical properties may differ.
We give criteria for such a tiling space to be weakly mixing, and for
the dynamics of two such spaces to be topologically conjugate. 
\end{abstract}

\maketitle

\markboth{Clark-Sadun}{When size matters: subshifts and tiling spaces}
\newpage

\section{Introduction}

We consider the dynamics of 1-dimensional minimal substitutions
subshifts (with a natural $\mathbb{Z}$ action) and the associated
1-dimensional tiling spaces (with natural $\mathbb{R}$ actions). Given
an alphabet $\mathcal{A}$ of $n$ symbols $\left\{
a_{1},...,a_{n}\right\} $, a \emph{substitution }on $ \mathcal{A}$ is
a function $\sigma $ from $\mathcal{A}$ into the non-empty, finite
words of $\mathcal{A}$. Associated with such a substitution is the $
n\times n$ matrix $M$ which has as its $\left( i,j\right) $ entry the
number of occurrences of $a_{j}$ in $\sigma \left( a_{i}\right) $. A
substitution is \emph{primitive }if some positive power of $M$ has
strictly positive entries. The substitution $\sigma $ induces the map
$\sigma :\mathcal{A}^{ \mathbb{Z}}\rightarrow
\mathcal{A}^{\mathbb{Z}}$ given by $\left\langle \dots
u_{-1}.u_{0}u_{1}\dots \right\rangle \overset{\sigma }{\mapsto }
\left\langle \dots \sigma \left( u_{-1}\right) .\sigma \left(
u_{0}\right) \sigma \left( u_{1}\right) \dots \right\rangle .$ For any
primitive substitution $\sigma $ there is at least one point $u\in
\mathcal{A}^{ \mathbb{Z}}$ which is periodic under $\sigma ,$ and the
closure of the orbit of any such $u$ under the left shift map $s$ of
$\mathcal{A}^{\mathbb{Z}}$ forms a minimal subshift
$\mathcal{S}$. This subshift $\mathcal{S}$ is uniquely determined by
$\sigma$.
To avoid trivialities, we shall only
consider primitive, \emph{aperiodic} substitutions $\sigma,$ 
i.e., those for which the
subshift $\mathcal{S}$ is not periodic.

Given a collection of intervals $\mathcal{I}=\left\{
I_{1},...,I_{n}\right\} $, a \emph{tiling} $T$ of $\mathbb{R}$ by
$\mathcal{I}$ is a collection of closed intervals $\left\{
T_{i}\right\} _{i\in \mathbb{Z}}$ satisfying

\begin{enumerate}
\item $\cup _{i\in \mathbb{Z}}T_{i}=\mathbb{R},$

\item For each $i\in \mathbb{Z},$ $T_{i}$ is the translate of some $I_{\tau
\left( i\right) }\in \mathcal{I},$ and

\item $T_{i}\cap T_{i+1}$ is a singleton for each $i\in \mathbb{Z}.$
\end{enumerate}

\noindent If the function $\tau :\mathbb{Z\rightarrow }\left\{
1,...,n\right\} $ is an element of a minimal substitution subshift
$\mathcal{ S}$ of $\left\{ 1,...,n\right\} ^{\mathbb{Z}}$, then $T$ is
called a \emph{substitution tiling.} There is a natural topology on
the space $\mathfrak{T}$ of tilings of $\mathbb{R}$ by $\mathcal{I}$ that is induced by
a metric which measures as close any two tilings $T$ and $T'$
that agree on a large neighborhood of $0$ up to an $\varepsilon $
translation. There is then the continuous translation action
$\mathbf{T}$ of $\mathbb{R}$ on $\mathfrak{T:}$ for $t\in \mathbb{R}$
and the tiling $T=\left\{ T_{i}\right\} _{i\in \mathbb{ Z}}$,
$\mathbf{T}:\left( T,t\right) \mapsto \mathbf{T}_{t}T=\left\{
T_{i}-t\right\} _{i\in \mathbb{Z}}.$ (A positive element of $\red$
moves the origin to the right, or equivalently moves tiles to the
left).  The closure of the translation orbit of any substitution
tiling $T$ in $\mathfrak{T}$ is then a minimal set of the action, the
(\emph{substitution})\emph{\ tiling space }$\mathcal{T}$ of the tiling
$T$. 
As the tiling space for any iterate of a substitution $\sigma $
is the same as the tiling space of $\sigma $ and since any
substitution has a point $u$ whose right half $\left\langle
u_{0}u_{1}\dots \right\rangle $ is periodic under substitution, for
ease of discussion we shall only consider substitutions with a point
$u$ whose right half $\left\langle u_{0}u_{1}\dots \right\rangle $ is
fixed under substitution.

Flows under a function provide an alternative description of tiling
spaces convenient for our purposes. Given a minimal subshift $\left(
\mathcal{S} ,s\right) $ of the shift on $\mathcal{A}^{\mathbb{Z}}$ and
$f:\mathcal{A} \rightarrow \left( 0,\infty \right) ,$ there is the
flow under $f$ given by the natural $\mathbb{R}$ action on
$\calT_{f}=\calS\times \mathbb{R}/\sim $, where $(u,f(u_{0}))\sim
(s(u),0)$. If we then associate with each $a_{i}\in \mathcal{A}$ a
closed interval $I_{i}$ of length $f\left( a_{i}\right) $ and form the
tiling on $T$ of $\mathbb{R}$ by $\left\{ I_{1},\dots ,I_{n}\right\} $
with associated function $\tau =u\in \mathcal{A}^{\mathbb{Z} }$, which
has the left endpoint of the interval corresponding to $u_{0}$ at $
0\in \mathbb{R}$, then the function sending $T$ to the class of
$\left( u,0\right) $ in $\calT_{f}$ extends uniquely to a
homeomorphism $\mathcal{T} \rightarrow \mathcal{T}_{f}$ which
conjugates the respective $\mathbb{R}$ actions.

The primary focus of this paper is the extent to which the dynamical
systems $\mathcal{T}_{f}$ depend on the function $f$. If all we care
about is the topological space, they don't:

\begin{thm}
\label{homeo}Let $\calS$ be a subshift, and let $f,g$ be two positive
functions on the alphabet of $\calS$. Then $\calT_{f}$ and $\calT_{g}$
are homeomorphic.
\end{thm}

\medskip \noindent \textbf{Proof.\enspace} Every class $x\in
\calT_{f}$ has a (unique) representative of the form $(u,t)$, with
$0\leq t<f(u_{0})$. Let $ h(x)$ be the class of
$(u,g(u_{0})t/f(u_{0}))$ in $\calT_{g}$. As the function on quotient
spaces induced by a continuous function $\mathcal{S} \times
\mathbb{R\rightarrow }\mathcal{S}\times \mathbb{R}$, $h$ is
continuous, and it has a well-defined continuous inverse.\hfill
$\square$

But what about dynamics? Dekking and Keane \cite{DK} showed that, for
the substitution $a\rightarrow abab$, $b\rightarrow bbba$, the
associated tiling space with $f(a)=2,f(b)=1$ has point spectrum only
at multiples of $2\pi $ . (That is, the associated discrete dynamical
system is weakly mixing.)  Further, Berend and Radin \cite{BR} showed
that a tiling space based on this substitution has trivial point
spectrum if $f(a)/f(b)$ is irrational. Radin and Sadun \cite{RS}
obtained very different results for the Fibonacci substitution
$a\rightarrow b$, $b\rightarrow ab$. For the Fibonacci substitution,
$\calT_{f}$ is topologically conjugate to $\calT_{g}$ whenever
$f(a)+\tau f(b)=g(a)+\tau g(b)$, where $\tau =(1+\sqrt{5})/2$ is the
golden mean. (The conjugacy is not the map $h$, described above, but
is homotopic to $h$. A generalization of this construction appears
below in Section 3.) In particular, all Fibonacci tiling spaces have
pure point spectrum, regardless of the sizes of the tiles.

In Section 2 we examine the general ergodic properties of the flows
$\calT _{f},$ showing that these flows are never mixing. We also
provide general criteria for the existence of point spectrum and
constraints on the form of that spectrum, in terms of the eigenvalues
of the substitution matrix and the ratios of lengths of the tiles. The
Dekking-Keane and Berend-Radin results are then easily understood. 

In Section 3 we provide criteria for $ \calT_{f}$ and $\calT_{g}$ to
be topologically conjugate.  In particular, if only one eigenvalue has
magnitude 1 or greater, we show that there exists a constant $c$
such that $\calT_f$ and $\calT_{cg}$ are topologically conjugate.
In addition to Pisot substitutions, this result applies to substitutions such
as Thue-Morse, whose matrices have zero as an eigenvalue.

Given a substitution, there are two natural choices of tile length.
One choice is the standard suspension, where all the tiles have length
one.  The other choice is to pick lengths according to the left
Perron-Frobenius eigenvector of the substitution matrix.  This yields
tilings that are self-similar in the sense of Solomyak
\cite{S}. 

For Pisot substitutions, there are long-standing conjectures that both
the self-similar tiling space and the suspended subshift have pure
point spectrum.  Partial results have
been obtained for both cases, using the ``balanced pair'' algorithm
\cite{BD2, HS} for the suspended subshift and the ``overlap algorithm''
for the self-similar tiling space \cite{S}.  Various relations
between the two problems have been derived \cite{SV}.  Our results 
imply that the two problems are in fact equivalent. 

The techniques of Section 3 draw heavily on the work of Moss\'e
\cite{M1, M2}. We generalize her computation of the number of times
a word can be repeated to allow for fractional repetitions, and show
that the set of period lengths yields a conjugacy invariant. The proof of
Theorem \ref{families} in particular is an extension of ideas in \cite{M1}. 
Moreover, we are implicitly using Moss\'e's bilateral recognizability
results whenever we refer to the location of supertiles within a tiling.
  
Finally, Sadun and Williams \cite{SW} have shown that every finite
dimensional tiling space meeting mild conditions (consisting of a
finite number of polygonal tiles, appearing in a finite number of
orientations, and meeting edge-to-edge) can be deformed into a tiling
space that is topologically conjugate to the $d$-fold suspension of a
${} \mathbb{Z}^{d}$ subshift. This deformation changes the shapes and
sizes of the tiles, but not the combinatorics of which tiles touch
which others.  The difference between the dynamics of general tiling
spaces and of subshifts is therefore not combinatorial in nature, but
depends rather on the geometry of the individual tiles.  This paper is
an exploration of this phenomenon in one dimension, where the geometry
of a tile is simply its length. We believe that the results of this
paper will provide essential tools for exploring the phenomenon in
higher dimensions as well.

\section{Ergodic Properties}

For a substitution tiling space $\mathcal{T}_{f},$ let
$L=(f(a_{1}),\ldots ,f(a_{n}))$ be the (row) vector that gives the
lengths of the various tiles.  Let $L_0$ be the normalized positive
left Perron-Frobenius eigenvector with eigenvalue $\lambda _{PF}$ of
the substitution matrix $M$, and let $\calT_{0}$ be the tiling space
associated with the length vector $L_{0}$. Every such tiling space
$\mathcal{T}_{0}$ is known to admit a self-homeomorphism derived from
$ \sigma $ (also denoted $\sigma )$ which inflates all the tiles by a
factor of $\lambda _{PF}$ and which maps the class of a point $\left(
x,0\right) $ to the class of $\left( \sigma \left( x\right) ,0\right)
,$ see, e.g., \cite {BD}. The more general substitution tiling space
$\mathcal{T}_{f}$ then admits a substitution homeomorphism, also
denoted $\sigma ,$ conjugate to the substitution homeomorphism of
$\mathcal{T}_{0}$ via the homeomorphism $h$ of Theorem \ref{homeo},
but this substitution homeomorphism is not generally affine. Solomyak
\cite{S} has investigated the ergodic properties of the self-affine
$\mathcal{T}_{0}$ in the case that the matrix $M$ is
diagonalizable. Now we proceed to investigate the ergodic properties
of the more general tiling spaces $\mathcal{T}_{f},$ requiring only
that the substitution be primitive and aperiodic. As most proofs in
this section are similar to those for analogous known results for
substitutions, the proofs shall be abbreviated.

It should be noted that most of the literature is devoted to the ergodic
properties of one-sided substitutions, while we are dealing with two-sided
substitutions. However, as there are finitely many $s-$orbits of words in
a one-sided substitution that admit more than one left extension \cite{Q}, 
(see \cite{BDH} for precise bounds), measure theoretic considerations are not 
altered. For example, by placing all two-sided words with the same right halves into the same partition
element, partitions of one-sided substitutions lead to partitions of
two-sided substitutions with corresponding partition elements having the
same measure.

\medskip \noindent \textbf{Definitions.\enspace }For a word $w=w_{0}\cdots
w_{k}$ from the alphabet $\mathcal{A}$ of the substitution subshift $
\mathcal{S}$ derived from the substitution $\sigma ,$ 
\begin{equation}
\left[ w\right] \overset{\text{def}}{=}\left\{ u\in \mathcal{S}
\,|\,u_{0}\cdots u_{k}=w\right\} .
\end{equation}
A \textit{cylinder set }$\left[ w\right] \times I$ in $\mathcal{T}_{f}$ for
a word $w=w_{0}\cdots w_{k}$ and interval $I\subseteq \left[ 0,f\left( w_{0}\right) \right) $ 
is the set of all $x\sim \left( u,t\right)
\in \mathcal{T}_{f}$ with $u\in \left[ w\right] $ and $t\in I.$

The \textit{population vector }$v=\left( v_{1},...,v_{n}\right) ^{T}$ of the
finite word $w$ gives the number of occurrences $v_{i}$ of the letter $a_{i}$
in $w.$  Note that $v$ is a column vector, not a row, and that the
substitution matrix acts from the left.

For $u\in \mathcal{A}^{\mathbb{Z}}$, a \textit{recurrence word }is a
finite word $w$ in $u,$ $w=u_{r}u_{r+1}\cdots u_{s}$ satisfying the
condition that $ u_{s+1}=u_{r}$, and a \textit{recurrence vector} is
the population vector for a recurrence word. For example, for $\ldots
abbcba\ldots $ in $\left\{ a,b,c\right\} ^{\mathbb{Z}}$, $abbcb$ is a
recurrence word for which $ (1,3,1)^{T}$ is the recurrence vector.
The vectors $(0,1,0)^{T}$, $(0,1,1)^{T}$, and $(0,2,1)^{T}$ are
recurrence vectors for the successive $b$'s.

A recurrence vector $v$ is called \textit{full} if the vectors
$\{M^kv\}$, with $k$ ranging from 0 to $n-1$, are linearly
independent.  If $M$ is diagonalizable, this is equivalent to $v$
pairing nontrivially with every left-eigenvector of $M$.  Note that
every primitive substitution on two letters admits a full recurrence vector,
namely $(0,1)^T$ or $(1,0)^T$.

Primitive aperiodic substitution tilings are recurrent, meaning that
any finite patch of one tiling appears somewhere in every other
tiling. As a result, the set of recurrence vectors is the same for
every tiling in the space, and so we can speak of the recurrence
vectors of a tiling \textit{space}.
\medskip
 
The following lemma and its implication for mixing are proved in much the
same way as for substitutions \cite[Theorem 2]{DK}.

\begin{lem} Let $v$ be the recurrence
vector of a recurrence word $r$ and let $t_{m}\overset{\text{def}}{=}
LM^{m}v$. For any cylinder set\textit{\ }$\left[ w\right] \times I$ in $
\mathcal{T}_{f}$ and for any $m$, let $S_{m}\overset{\text{def}}{=}\left( 
\left[ w\right] \times I\right) \cap \left( \mathbf{T}_{-t_{m}}\left( \left[
w\right] \times I\right) \right) .$ Then 
\begin{equation}
\underset{m\rightarrow \infty }{\lim \inf }\;\mu \left( S_{m}\right) \geq
\gamma\cdot \mu \left( \left[ w\right] \times I\right) ,
\end{equation}
where $\gamma$ is a positive constant independent of $w$ and $I,$ depending
only on the recurrence word $r$ of $v.$
\label{olap} \end{lem}

\demo{Proof.}  With $t_{m}'\overset{ \text{def}}{=}\left(
1,\dots ,1\right) M^{m}v,$ it follows from \cite{DK} that the measures
$\mu'$ in the substitution subshift of the sets $
S_{m}'\overset{\text{def}}{=}\left[ w\right] \cap
s^{-t_{m}'}\left[ w\right] $ satisfy
\begin{equation}
\underset{m\rightarrow \infty }{\lim \inf }\;\mu'\left(
S_{m}'\right) \geq \gamma'\cdot \mu'\left( 
\left[ w\right] \right) ,
\end{equation}
for some positive constant $\gamma'$ independent of $w.$
Observe that $\mu'\left( \left[ w\right] \right) $ is the
limit as $ k\rightarrow \infty $ of the average number of iterates
$\left\{ s^{i}\left( u\right) |\,i=0,...,k-1\right\} $ in $\left[
w\right] $ $,$ for any given $ u\in \mathcal{A}^{\mathbb{Z}}$. And the
measure $\mu \left( \left[ w\right] \times I\right) $ in
$\mathcal{T}_{f}$ is the limit as $T\rightarrow \infty $ of the
average amount of time that $\left\{ \mathbf{T}^{t}\left( x\right)
|\,t\in \lbrack 0,T]\right\} $ spends in $\mu \left( \left[ w\right]
\times I\right) $ for any given $x\in \mathcal{T}_{f}.$ Hence, $\mu
\left( \left[ w \right] \times I\right) =\mu'\left( \left[
w\right] \right) \times \dfrac{\text{length}\left( I\right) }{\mu
'\left( \left[ a_{1} \right] \right) f\left( a_{1}\right)
+\cdots +\mu'\left( \left[ a_{n}\right] \right) f\left(
a_{n}\right) }.$ Now, if $u\in \left[ w\right] $ and $t\in I,$ then
$s^{t_{m}'}\left( u\right) \in \left[ w\right] $ if and only
if $\mathbf{T}_{t_{m}}\left( \left( u,t\right) \right) \in \left[
w\right] \times \left\{ t\right\}$, and the result follows.

\hfill $\square$\medskip

By choosing the cylinder set smaller in measure than a $\gamma$ for a
fixed recurrence word, this leads to the following conclusion.

\begin{thm} None of the substitution tiling space
flows on $\mathcal{T}_{f}$ are (strongly) mixing.
\end{thm}

\begin{thm}

The number $k$ is in the point spectrum of $\mathcal{T}_{f}$ if and
only if, for every recurrence vector $v$,
\begin{equation}
\dfrac{k}{2\pi }LM^{m}v\rightarrow 0 \text{(mod 1) as } m\rightarrow
\infty .
\end{equation}
Moreover, when this condition is met there is a constant $C$ and a $\rho <1$
satisfying
\begin{equation}
\left| \exp \left( ikLM^{m}v\right) -1\right| <C\rho ^{m}
\end{equation}
for any given $v.$

\label{zeromod1}\ \end{thm} \medskip

\demo{Proof.}

Essentially the same arguments as found in \cite{FMN} using Rokhlin
stacks and columns apply to show that the condition is both necessary
and sufficient for \emph{return words, }finite words $w$ in $u$, $
w=u_{r}u_{r+1}\cdots u_{s}$ satisfying for sufficiently large $m$ the
conditions $\sigma ^{m}u_{r}=\sigma ^{m}u_{s+1}$ and $\sigma
^{m}u_{r}\neq \sigma ^{m}u_{j}$ for $j=r+1,...,s$. Since the
eigenfunction constructed as in \cite{FMN} using the geometric
convergence and the condition for return words is continuous, one need
only consider continuous eigenfunctions.

We now show that for continuous eigenfunctions the condition for the
general \emph{recurrence }word is necessary, which is all that need be
shown since the condition with return words is already
sufficient. Consider then a continuous eigenfunction $f$ with
eigenvalue $k$. Let $x\sim \left( u,0\right) \in \mathcal{T}_{f}$,
where $u\in \mathcal{A}^{\mathbb{Z}}$ is \ a point of the substitution
$\sigma $ with fixed right half $\left\langle u_{0}u_{1}\dots
\right\rangle $, and let $v$ be a recurrence vector for a recurrence
word between $u_{0}$ and some $u_{s}$, and let $t_{m}\overset{
\text{def}}{=}LM^{m}v$. Note that $\mathbf{T}_{t_{0}}x$ agrees with
$x$ on the interval covered by the single tile of $x$ represented by
$u_{0}$. Then $ \mathbf{T}_{t_{m}}x=\sigma ^{m}(\mathbf{T}_{t_{0}}x)$
agrees with $x$ on forward time intervals of increasing unbounded
length, and so $\exp (ikt_{m}) $ must converge to 1, and $kt_{m}/(2\pi
)$ must converge to 0 (mod 1). Now let $v_{1}$ and $v_{2}$ be two such
recurrence vectors. Since $ kLM^{m}v_{1}/(2\pi )$ and
$kLM^{m}v_{2}/(2\pi )$ both converge to zero (mod 1), so does
$kLM^{m}(v_{1}-v_{2})/(2\pi )$. In other words, we can use any
recurrence vector for the type of letter (say, $a$) that occurs at
$u_{0}$.  Finally, suppose that $v_{3}$ is a recurrence vector for
letter $b$. Since $ \sigma $ is primitive, suppose $a$ appears
somewhere in the $\ell $-th substitution of $b$, for some integer
$\ell $. Then $v_{4}=M^{\ell }v_{3}$ is a recurrence vector for $a$,
so $kLM^{m}v_{3}/(2\pi )=kLM^{m-\ell }v_{4}/(2\pi )$ converges to
$0($mod $1)$. To obtain the geometric convergence as in the statement
of the theorem, one applies the argument of \cite[Lemme 1]{H} with the
lattice $f\left( a_{1}\right) \mathbb{Z\oplus \cdots }\mathbb{\oplus }
f\left( a_{n}\right) \mathbb{Z}$ in place of the lattice
$\mathbb{Z}^{n}.$ \hfill $\square$\medskip

The application of this criterion depends on the eigenvalues of $M$
and on the possible forms of the recurrence vectors.

\begin{thm}
Suppose that all the eigenvalues of $M$ are of magnitude 1 or greater,
and that there exists a full recurrence vector. If the ratio of any
two tile lengths is irrational, then there is trivial point
spectrum. If the ratios of tile lengths are all rational, then the
point spectrum is contained in $2\pi\qed/L_{1}$.
\label{allbig1}
\end{thm}

\medskip \noindent \textbf{Proof.\enspace} Let $k$ be in the point
spectrum, and consider the sequence of real numbers
$t_{m}=kLM^{m}v/(2\pi )$, where $v$ is a fixed full recurrence
vector. Let $p(\lambda )=\lambda ^{n}+a_{n-1}\lambda ^{n-1}+\cdots
+a_{0}$ be the characteristic polynomial of $M$. Note that the
$a_{i}$'s are all integers, since $M$ is an integer matrix. Since
$p(M)=0$, the $t_{m}$'s satisfy a recurrence relation:
\begin{equation}
t_{m+n}=-\sum_{k=0}^{n-1}a_{k}t_{m+k}.  \label{intrecur}
\end{equation}
By Theorem \ref{zeromod1}, the $t_{m}$'s converge to zero $\pmod{1}$. 
That is, we can write 
\begin{equation}
t_{m}=i_{m}+\epsilon _{m}  \label{splits}
\end{equation}
where the $i_{m}$'s are integers, and the $\epsilon _{m}$'s converge
to zero as real numbers. By substituting the division (\ref{splits})
into the recursion (\ref{intrecur}), we see that both the $i$'s and
the $\epsilon $'s must separately satisfy the recursion
(\ref{intrecur}), once $m$ is sufficiently large (e.g., large enough
that the $\epsilon $'s are bounded by $1/\sum |a_{i}|$). However, any
solution to this recursion relation is a linear combination of powers
of the eigenvalues of $M$ (or polynomials in $m$ times eigenvalues to
the $m$-th power, if $M$ is not diagonalizable).  Since the
eigenvalues are all of magnitude one or greater, such a linear
combination converges to zero only if it is identically
zero. Therefore $ \epsilon _{m}$ must be identically zero for all
sufficiently large values of $m$.

Each $t_{m}$ is an integer linear combination of the elements of the
vector $ kL/(2\pi )$. From a sequence of $n$ consecutive $t$'s, one
can recover all the elements of $kL/(2\pi)$ by inverting an integer matrix
whose columns are powers of $M$ times $v$. (The invertibility of this
matrix depends on the fact that $v$ is full and that zero is not an
eigenvalue.) By Cramer's rule, this inversion involves integer
multiplication and addition, followed by division by the determinant
of this matrix. Since the $t_{m}$'s are integers (for $m$ large
enough), the components of $kL/(2\pi )$ must then all be rational. If
$L_{i}/L_{j}$ is irrational, this implies that $k=0$. If all the $L$'s
are rationally related, it implies that $k$ lies in $2\pi\qed/L_{1}$
.\hfill $\square$\medskip

Theorem \ref{allbig1} gives partial spectral information. When the alphabet
has only two letters and the substitution has constant length, we can say
considerably more:

\begin{thm}
\label{constantlength} Suppose that we have a primitive substitution on two
letters of constant length $n$, where $\sigma(a)$ contains $n_a$ $a$'s
and $ n-n_a$ $b$'s, while $\sigma(b)$ contains $n_b$ $a$'s and $n-n_b$
$b$'s.  Suppose further that $1 \le n_a,n_b \le n-1$ and $n_a \ne
n_b$. Let $z$ be the greatest common factor of $n$ and $n_a-n_b$. Then
the point spectrum $ \sigma_{pp}$ depends as follows on the ratio
$L_1/L_2$:

\begin{enumerate}
\item If $L_{1}=L_{2}$, then there is a positive integer $N$ such that 
\begin{equation}
N{}\mathbb{Z}[1/n]\subset N{}L_{1}\sigma _{pp}/2\pi \subset \mathbb{Z}
{}[1/n].
\end{equation}

\item If $L_{1}/L_{2}\in \mathbb{Q}-\{1\}$, then there exist positive
integers $N_{1}$ and $N_{2}$ such that 
\begin{equation}
N_{1}{}\mathbb{Z}[1/z]\subset N_{2}{}L_{1}{}\sigma _{pp}/2\pi \subset {}
\mathbb{Z}[1/z].
\end{equation}

\item If $L_{1}/L_{2}\notin \mathbb{Q}$, then $\sigma _{pp}=\{0\}$.
\end{enumerate}
\end{thm}

\medskip \noindent \textbf{Proof.\enspace} The $L_{1}=L_{2}$ case is
not new -- Coven and Keane \cite{CK} proved an even stronger result
over 30 years ago. We prove all three cases together, since the proof
neatly illustrates the significance of the ratio $L_{1}/L_{2}$.

The substitution matrix 
\begin{equation}
M = 
\begin{pmatrix}
n_a & n_b \cr n-n_a & n-n_b
\end{pmatrix}
\end{equation}
has eigenvalues $n$ and $n_a-n_b$, with left-eigenvectors $(1,1)$ and
$ (n-n_a, -n_b)$ and right-eigenvectors $(n_b, n-n_a)^T$ and $(1,
-1)^T$.  Since $n_a \ne n_b$, both eigenvalues have magnitude one or
greater.

Since the system is aperiodic and there are only two letters, either $
(1,0)^{T}$ or $(0,1)^{T}$ is a (full) recurrence vector. Suppose that
$(1,0)^{T}$ is (the other case is similar). Then we take

\begin{equation}
v = 
\begin{pmatrix}
1 \\ 
0
\end{pmatrix}
= \frac{1}{n+n_b-n_a} 
\begin{pmatrix}
n_b \\ 
n-n_a
\end{pmatrix}
+ \frac{n-n_a}{n+n_b-n_a} 
\begin{pmatrix}
1 \\ 
-1
\end{pmatrix}
.
\end{equation}

As in the proof of Theorem \ref{allbig1}, if $k\in \sigma _{pp}$ then the
numbers $t_{m}=kLM^{m}v/2\pi $ must eventually be integers. Note that 
\begin{equation}
\begin{split}
2\pi t_{m}=& k 
\begin{matrix}
LM^{m}
\end{matrix}
\begin{pmatrix}
1 \\ 
0
\end{pmatrix}
\\
2\pi t_{m+1}=& k 
\begin{matrix}
LM^{m+1}
\end{matrix}
\begin{pmatrix}
1 \\ 
0
\end{pmatrix}
=k 
\begin{matrix}
LM^{m}
\end{matrix}
\begin{pmatrix}
n_{a} \\ 
n-n_{a}
\end{pmatrix}
,
\end{split}
\end{equation}
so 
\begin{equation}
2\pi (t_{m},t_{m+1})=kLM^{m} 
\begin{pmatrix}
1 & n_{a} \\ 
0 & n-n_{a}
\end{pmatrix}
,
\end{equation}
and 
\begin{equation}
kLM^{m}=2\pi (t_{m},t_{m+1}){\ 
\begin{pmatrix}
1 & n_{a} \\ 
0 & n-n_{a}
\end{pmatrix}
}^{-1}=\frac{2\pi}{n-n_{a}}\left( (n-n_{a})t_{m},t_{m+1}-n_{a}t_{m}\right) .
\end{equation}
Multiplying on the right by the right-eigenvectors of $M$ gives two scalar
equations for $k$: 
\begin{equation}
\begin{split}
k(L_{1}-L_{2})(n_{a}-n_{b})^{m}=& \frac{2\pi }{n-n_{a}}(nt_{m}-t_{m+1}) \\
k(n_{b}L_{1}+(n-n_{a})L_{2})n^{m}=& 2\pi ((n_{b}-n_{a})t_{m}+t_{m+1}).
\end{split}
\label{kequations}
\end{equation}

If $L_{1}=L_{2}$, the first equation reads $0=0$, but the second equation
says that $NkL_{1}n^{m}/2\pi $ is an integer for $m$ sufficiently large,
where $N=n+n_{b}-n_{a}$. Thus $NL_{1}\sigma _{pp}/2\pi \subset \mathbb{Z}
[1/n]$.

If $L_{1}/L_{2}$ is rational but not one, then the first equation says
that $ \sigma _{pp}L_{1}/2\pi $ is contained in a fixed rational
multiple of $ \mathbb{Z}[1/(n_{a}-n_{b})]$, while the second says that
$\sigma _{pp}L_{1}/2\pi $ is contained in a fixed rational multiple of
$\mathbb{Z} [1/n]$. Thus an integer multiple of $\sigma
_{pp}L_{1}/2\pi $ lies in ${} \mathbb{Z}[1/n]\cap
\mathbb{Z}[1/(n_{a}-n_{b})]=\mathbb{Z}[1/z]$.

If $L_{1}/L_{2}$ is irrational, then the only simultaneous solution to
both equations is $k=0$. For if $k\neq 0$ is a solution, then the
ratio of the two equations would imply that
$(L_{1}-L_{2})/(n_{b}L_{1}+(n-n_{a})L_{2})$ would be rational, which
is inconsistent with the irrationality of $ L_{1}/L_{2}$. Thus $\sigma
_{pp}\subset \{0\}$ and so $\sigma _{pp}=\{0\}.$ All that remains is
to show that $L_{1}\sigma _{pp}/2\pi $ contains ${} \mathbb{Z}[1/n]$
if $L_{1}=L_{2}$, and that $L_{1}\sigma _{pp}/2\pi $ contains a
rational multiple of ${}\mathbb{Z}[1/z]$ if $L_{1}/L_{2}$ is rational
but not one.

We call the result of applying the substitution $\ell $ times to a
single tile (or a single letter) a \textit{supertile} of order $\ell
$. If $ L_{1}=L_{2}$, then all supertiles of order $m$ have length
$n^{m}L_{1}$.  Thus the coordinates of the endpoints of such
supertiles agree modulo $ n^{m}L_{1}$. For each pair of integers
$(j,m)$ we define a function on our tiling space. If $x$ is a tiling,
let
\begin{equation}
\psi _{j,m}(x)=\exp (2\pi ijp/n^{m}L_{1}),
\end{equation}
where $p$ is the endpoint of \textit{any} supertile of order $m$ in $x$.
This is manifestly an eigenfunction of translation with eigenvalue 
$2 \pi j/ n^{m} L_1$.

If $L_{1}/L_{2}$ is rational but not integral, then the length of each
supertile of order $m$ is a fixed rational linear combination of $
L_{1}n^{m} $ and $L_{1}(n_{a}-n_{b})^{m}$. Thus there is a constant
$c$, a rational multiple of $L_{1}$, such that every supertile of
order $m$ has length divisible by $cz^{m}$. For each pair $(j,m)$ we
then define the eigenfunction
\begin{equation}
\phi _{j,m}(x)=\exp (2\pi ijp/cz^{m}),
\end{equation}
which shows that some rational multiple of $\mathbb{Z}[1/z]$ is
contained in $L_{1}\sigma _{pp}/2\pi $.\hfill $\square$

\medskip

As a special case, we recover the results of Dekking and Keane, and
then of Berend and Radin, since their substitution $\sigma $ has $
n=4,n_{a}=2,n_{b}=1 $, so $z=1$. When $L_{1}=2$ and $L_{2}=1$, the
point spectrum is $2\pi \mathbb{Z}[1/z]=2\pi \mathbb{Z}$. This also
illustrates the sharp contrast between substitution tiling spaces
admitting a flow with point spectrum and spaces admitting
equicontinuous flows. If $X$ admits an equicontinuous $\mathbb{R}$
action $\phi $ with point spectrum $\sigma _{pp}$ and if $\phi
'$ is any $\mathbb{R}$ action on $X$ obtained by a time
change of $\phi ,$ then the point spectrum of $\phi '$ is
$c\cdot \sigma _{pp}$ for some possibly $0$ constant $c,$ see Egawa
\cite {E}. In the substitution $\sigma $ currently under
consideration, the $ \left( 1,1\right) $ flow has purely point
spectrum $\sigma _{pp}$ but admits a time change represented by
$\left( 2,1\right) $ which has discrete spectrum which can only be
represented as a proper subset of a multiple of $ \sigma _{pp}.$

What about the Fibonacci tiling? In that case the smaller eigenvalue is $
1-\tau$, whose magnitude is less than 1, so Theorem \ref{allbig1} does not
apply directly. However, there is still something that can be said about
such cases:

\begin{thm}
\label{somesmall1} Suppose that there exists a full recurrence vector. 
Let $S$ be the span of the (generalized left-) eigenspaces of $M$ with
eigenvalue of magnitude strictly less than 1. If $k$ is in the point
spectrum, then $kL/ {2 \pi}$ is the sum of a rational vector and an
element of $S$.
\end{thm}

\medskip\noindent\textbf{Proof.\enspace} 
Let $v$ be a full recurrence vector, and let $t_m = kLM^mv/2\pi$, as 
in the proof of Theorem \ref{allbig1}. As before, $t_m$
converges to zero (mod 1), so we can write $t_m = i_m + \epsilon_m$
with $i_m$ integral and $\epsilon_m$ converging to zero, and with both
the $i_m$'s and the $\epsilon_m$'s eventually satisfying the recursion
relation (\ref {intrecur}). Since the $\epsilon_m$'s converge to zero,
they must eventually be linear combinations of $m$-th powers of the
small eigenvalues of $M$ (possibly times polynomials in $m$, if $M$ is
not diagonalizable).  By adding an element of $S$ to $L$ (and hence to
$kL/2\pi$), we can then get all the $\epsilon_m$'s to be identically
zero beyond a certain point.

If $M$ is invertible, the resulting value of $kL/2\pi$ must then be
rational by the same argument as in the proof of Theorem
\ref{allbig1}. If $M$ is not invertible, then we need an extra step,
since the vectors $M^m v, M^{m+1} v, \ldots, M^{m+n-1}v$ are no longer
linearly independent.  Let $S_0 \subset S$ be the generalized
left-eigenspace with eigenvalue zero, and suppose that this space has
dimension $\ell$. The vectors $M^mv, \ldots, M^{m+n-\ell-1}v$ are
linearly independent and span the annihilator of $S_0$, which can be
identified with the dual space of ${\red^n}^*/S_0$.  From the integers
$t_m, \ldots, t_{m+n-\ell-1}$ we can reconstruct, by rational
operations, a representative of $kL/2\pi$ in ${\red^n}^*/S_0$.  But
that implies that $kL/2\pi$ is a rational vector plus an element of
$S_0$.  \hfill $\square$
\medskip

Another way of stating the same result is to say that $kL/2\pi$,
projected onto the span of the large eigenvectors, equals the
projection of a rational vector onto this span.

This theorem can be used in two different ways. First, it constrains
the set of vectors $L$ for which the system admits point spectrum. Let
$d_b$ be the number of large eigenvalues, counted with (algebraic)
multiplicity.  There are only a countable number of possible values
for the projection of $kL/2\pi$ onto the span of the large
(generalized) eigenvectors.  If the projection of $L$ does not lie in
the ray generated by one of these points, there is no point
spectrum. In other words, one must tune $d_b-1$ parameters to a
countable number of possible values in order to achieve a nontrivial
point spectrum.

A second usage is to constrain the spectrum for fixed $L$.  The
rational points in $\red^n$, projected onto the span of the large
eigenvalues, and then intersected with the ray defined by a fixed $L$,
form a vector space over $\qed$ of dimension at most $n+1-d_b$. As a
result, the point spectrum tensored with $\mathbb{Q}$ is a vector
space over $\qed$ whose dimension is bounded by one plus the number of
small eigenvalues.  Below we derive an even stronger result, in which
only the small eigenvalues that are conjugate to the Perron-Frobenius
eigenvalue contribute to the complexity of the spectrum.

So far we have assumed that there exists a full recurrence vector $v$.
This may fail, either because there is something peculiar about the
available recurrence vectors, or because the matrix $M$ has repeated
eigenvalues.  In such cases we can repeat our analysis on a proper
subspace of ${}\mathbb{R}^{n}$. There are necessarily fewer
constraints on the vector $L$, but for fixed $L$ the constraints on
$k$ are actually stronger than before.

\begin{thm}
\label{notfull} Let $b_{PF}$ be the number of large eigenvalues that
are algebraically conjugate to the Perron-Frobenius eigenvalue
$\lambda_{PF}$ (including
$\lambda_{PF}$ itself), and let
$s_{PF}$ be the number of small eigenvalues conjugate to $\lambda_{PF}$. 
For the system to have nontrivial point spectrum, $L$ must lie in a 
countable union of subspaces of $\red^n$, each of codimension $b_{PF}-1$
or greater.  For fixed $L$, the dimension over $\qed$ of the point
spectrum tensored with $\qed$ is at most $s_{PF}+1$.  
\end{thm}

\medskip\noindent\textbf{Proof.\enspace} As a first step we
diagonalize $M$ over the rationals as far as possible.  By rational
operations we can always put $M$ in block-diagonal form, where the
characteristic polynomial of each block is a power of an irreducible
polynomial.  Since the Perron eigenvalue $\lambda _{PF}$ has algebraic
multiplicity one, every eigenvalue algebraically conjugate to $\lambda
_{PF}$ also has multiplicity one. Thus there is a unique block whose
characteristic polynomial has $\lambda _{PF}$ for a root. Every
recurrence vector is integral, and pairs nontrivially with $L_0$, and
so pairs nontrivially with every left-eigenvector of this block.  We
then consider the constraints on the spectrum that can be obtained
from this block alone.

We repeat the argument of the previous two theorems.  As before, we
modify $L$ by adding a linear combination of small eigenvectors so as
to make $t_m$ eventually integral.  Since the block diagonalization
was over the rationals, the vector $\vec t= (t_m, \ldots,
t_{m+n-1})^T$, expressed in the new basis, is still rational. From
this vector we deduce, as before, the projection of $kL/(2\pi)$ on any
large eigenvector with which our recurrence vector $v$ pairs
nontrivially.  This gives at least $b_{PF}$ independent constraints on
the pair ($L,k$), hence $b_{PF}-1$ constraints on $L$ if the spectrum
is to be nontrivial.

Next we consider the projection of $kL/(2\pi)$ onto the large
eigenspaces of the Perron-Frobenius block.  Since only $b_{PF}+s_{PF}$
components of $\vec t$ (expressed in the new basis) contribute, this
is the projection of $\qed^{b_{PF}+s_{PF}}$ onto $\red^{b_{PF}}$,
whose real span is all of $\red^{b_{PF}}$.  Intersected with the ray
defined by a fixed $L$, this gives a vector space of dimension at most
$s_{PF}+1$ in which $k$ can live.  \hfill $\square$

\section{Topological conjugacies}

In Section 2 we studied the point spectrum of a tiling space. This is
a conjugacy invariant; systems with different discrete spectra cannot
be conjugate. However, we found many different values of $L$ that gave
exactly the same point spectrum -- namely $\{0\}$. In this section we
give criteria for determining when two values of $L$ give rise to
conjugate tiling spaces. We supply sufficient conditions for conjugacy
in terms of the matrix $M$. We also supply necessary conditions; these
require knowing something about the actual substitution, and not just
its matrix $M$.  Finally, we exhibit examples that demonstrate the
need for such details.  Different substitutions with the same matrix
may result in different criteria for conjugacy.

Throughout this section, $\calS$ is a substitution subshift with
substitution map $\sigma $ and substitution matrix $M$ and we consider
two possible length vectors $L=(f(a_{1}),f(a_{2}),\ldots )$ and
$L'=(g(a_{1}),g(a_{2}),\ldots )$.  The substitution applied
$\ell $ times to a single tile (or letter) is called a
\textit{supertile} of order $\ell$.

\begin{thm}
If, for some integer $k$, 
\begin{equation}
\lim_{m \to \infty} \left ( L M^{m+k} - L^{\prime}M^m \right ) = 0,
\label{samesize}
\end{equation}
then $\calT_f$ is conjugate to $\calT_g$. \label{sufficient}
\end{thm}

\begin{cor}
If a substitution has only one eigenvalue of magnitude 1 or greater
(counted with multiplicity), then for arbitrary $f$ and $g$ there
exists a constant $c$ such that $\calT_f$ is conjugate to
$\calT_{cg}$. \label{Pisot}
\end{cor}

\medskip\noindent\textbf{Proof of corollary.\enspace} Recall that
$L_0$ is the normalized left Perron-Frobenius eigenvector of $M$. By
rescaling $g$ if necessary, we can assume that $L-L'$ is a linear
combination of (generalized) left eigenvectors other than $L_0$.
However, all the remaining eigenvectors have eigenvalue less than one,
so condition (\ref{samesize}) is met with $k=0$.  \hfill$\square$

\medskip\noindent\textbf{Proof of Theorem.\enspace} We prove the
theorem in two steps.  First we prove it in the special case that
$k=0$, and then in the case that $ L = L^{\prime}M$. Combining these
results then gives the theorem for all $k \le 0$. The situation for $k
>0$ is the same, only with $\calT_f$ and $\calT _g$ reversed.

If (\ref{samesize}) holds with $k=0$, then the supertiles in the
$\calT_{f}$ system asymptotically have the same size as the
corresponding supertiles in the $\calT_{g}$ system, and the
convergence is exponential. We then adapt the argument that Radin and
Sadun \cite{RS} applied to the Fibonacci tiling.  If $x$ is a tiling
in $\calT_{f}$ with tiles in a sequence $u=\ldots
u_{-1}u_{0}u_{1},\ldots $, then we require $\phi _{\ell }(x)\in
\calT_{g}$ to be a tiling with the exact same sequence of tiles. The
only question is where to place the origin. Let $d_{\ell }$ be the
coordinate of the right edge of the order-$\ell $ supertile of $x$
that contains the origin, and let $-e_{\ell }$ be the coordinate of
the left edge. Place the origin in $\phi _{\ell }(x)$ a fraction
$e_{\ell }/(d_{\ell }+e_{\ell })$ of the way across the corresponding
supertile of order $\ell $ in the $\calT_{g}$ system.  (This is
essentially the map $h$ of Theorem \ref{homeo}, only applied to
supertiles of order $\ell $.) Since the sizes of the supertiles
converge exponentially, the location of the origin in $\phi _{\ell
}(x)$ converges, and we can define $\phi (x)=\lim_{\ell \rightarrow
\infty }\phi _{\ell }(x)$. All that remains is to show that $\phi $
is a conjugacy. The approximate map $\phi _{\ell }$ is not a
conjugacy; translations in $x$ that keep the origin in the same
supertile of order $\ell $ get magnified (in $\phi (x)$) by the ratio
of the sizes of the supertiles of order $\ell $ containing the origin
in $x$ and $\phi (x)$. However, this ratio goes to one as $\ell
\rightarrow \infty $, while the range of translations to which this
ratio applies grows exponentially. In the $\ell \rightarrow \infty $
limit, $\phi $ commutes with translation by any $s\in (-e_{\infty
},d_{\infty })$.  Typically this is everything. The case where $x$
contains two infinite-order supertiles is only slightly trickier --
once one sees that $\phi $ preserves the boundary between these
infinite-order supertiles, it is clear that $\phi $ commutes with all
translations.

Now suppose that $L=L'M$. Then the tiles in the $\calT_{f}$
system have exactly the same size as the corresponding supertiles of
order 1 in the $\calT_{g}$ system. We define our conjugacy $\phi $ as
follows: If $x$ is a tiling in $\calT_{f}$ with tiles in a sequence
$u=\ldots u_{-1}u_{0}u_{1}\ldots $, then $\phi (x)$ to a tiling in
$\calT_{g}$ with sequence $\sigma (u)$, with each tile in $x$ aligned
with the corresponding order-1 supertile in $\phi (x)$.

Combining the cases, we obtain conjugacy whenever (\ref{samesize})
applies.\hfill $\square $\medskip

Theorem \ref{sufficient} gives sufficient conditions for two tiling
spaces to be conjugate. To derive necessary conditions we must
understand the extent to which words in our subshift $\calS$ (and its
associated tiling spaces) can repeat themselves. We will show that
words that repeat themselves $p>1$ times or more can be grouped into a
finite number of families, and we construct a conjugacy invariant from
the asymptotic lengths of these words.  Comparing these invariants for
different tile lengths then gives necessary conditions for conjugacy
(Theorem \ref{necessary}, below).

Moss\'{e} shows in
\cite{M1} that in every primitive aperiodic substitution subshift
there exists an integer $N$ so that no word is ever repeated $N$ or
more times. We will need to refine these results by considering words
that are repeated a fractional number of times. Counting fractional
degrees depends on the length vector: the word $ababa$ has its basic
period $ab$ repeated $2+{f(a) \over f(a)+f(b)}$ times.

\noindent \textbf{Definition.} A recurrence vector $v$ is a \textit{
repetition vector of degree }$p$\textit{\ }of a tiling space $\calT$ if
there are finite words $w$ (with population vector $v$) and $w'$,
such that: 1) $w'$ contains $w$, 2) $w'$ is periodic with
period equal to the length of $w$, 3) the length $Lw'$ of $w'$ is at
least $p$ times the length $Lw$ of $w$, and 4) $w'$ appears in some (and
therefore every) tiling in $\calT$.\smallskip

Note that every repetition vector of degree $p>p^{\prime}$ is also a
repetition vector of degree $p^{\prime}$. Note also that the word $w$
is typically not uniquely defined by $v$. A cyclic permutation of the
tiles in $w$ typically yields a word that works as well. For
instance, if $ w^{\prime}=ababa$ appears in a tiling, then $v=(1,1)^T$
is a repetition vector with degree $2+ {f(a) \over f(a)+f(b)}$, and we may
take either $w=ab$ or $ w=ba$.

Recall that  $\calT_0$ is the tiling space with length vector $L_0$, where
$L_0$ is the left Perron-Frobenius eigenvector of M.
For any word $w$, we denote by $|w|$ the length of the word $w$ in 
$\calT_{0}$. That is, $|w| = L_0 v$, where $v$ is the population vector
of $w$. If $v$ is a repetition vector of
degree $p$ for $\calT_{0}$, representing the word $w$ sitting inside $
w'$ in $u\in \mathcal{A}^{\mathbb{Z}}$, then $Mv$ is a repetition
vector of degree $p$, representing the word $\sigma (w)$ sitting inside $
\sigma (w')$ in $\sigma (u)\in \mathcal{A}^{\mathbb{Z}}$. The
degree is the same, since in $\calT_{0}$ the substitution $\sigma $
stretches each word by exactly the same factor, namely $\lambda _{PF}$.
Thus, every repetition vector $v$ gives rise to an infinite family of
repetition vectors $M^{k}v$. The following theorem limits the number of such
families.

\begin{thm}
Let $p>1$. There is a finite collection of vectors $\left\{ v_{1},\ldots
,v_{N}\right\} $ such that every repetition vector of degree $p$ for $\calT
_{0}$ is of the form $M^{k}v_{i}$ for some pair $(k,i)$. \label{families}
\end{thm}

\medskip \noindent \textbf{Proof.\enspace} From Moss\'e \cite{M1, M2},
we know that there is a recognition length $D_{ \mathcal{S}}$ of the
subshift $\calS$ such that knowing a letter, its $D_{_{ \mathcal{S}}}$
immediate predecessors and its $D_{_{\mathcal{S}}}$ immediate
successors determines the supertile of order 1 containing that letter,
and the position within that supertile of the letter. Thus for every
substitution tiling space $\calT$ there is a recognition length
$D_{\mathcal{ T}}$ such that the neighborhood of radius
$D_{_{\mathcal{T}}}$ about a point determines the supertile of order 1
containing that point, and the position of that point within the
supertile. (E.g., one can take $D_{_{\mathcal{T}}}$ to be
$1+D_{_{\mathcal{S}}}$ times the size of the largest tile.) Let $D=D_{
\mathcal{T}_{0}}$.

Suppose $v$ is a repetition vector of degree $p$ for $\mathcal{T}_{0}$,
corresponding to a word $w$ that is repeated $p$ times in a word $w'$, 
whose length is much greater than $D$. Then there is an interval of
size $|w'|-2D$ in which the supertiles of order 1 are periodic with
period $|w|$. A cyclic permutation of $w$ is therefore the union of 
distinct supertiles of order 1. Thus there exists a population vector 
$v_{1}$ of a word 
$w_{1}$, such that $v=Mv_{1}$, and such that $v_{1}$ is a repetition vector
of degree $(|w'|-2D)/|w|=p-2D/\lambda _{PF}|w_{1}|$. Note that $w$
itself does not have to equal $\sigma (w_{1})$ --- it may happen that $w$ is
a cyclic permutation of $\sigma (w_{1})$.

Repeating the process, we find words $w_i$ such that $\sigma(w_i)$ equals $
w_{i-1}$ (up to cyclic permutation of tiles), and such that the population
vector of $w_i$ is a repetition vector of degree 
\begin{equation}
p - \frac{2 D}{|w_i|} (\lambda_{PF}^{-1} + \lambda_{PF}^{-2} + \cdots +
\lambda_{PF}^{-i}) \ge p - \frac{2 D}{|w_i|} \sum_{\ell=1}^\infty
\lambda_{PF}^{-\ell} = p - \frac{2D}{|w_i| (\lambda_{PF}-1)}.
\end{equation}

Now pick $\epsilon <p-1$. We have shown that every repetition vector of
degree $p$ is of the form $M^{k}v_{i}$, where $v_{i}$ is a repetition vector
of degree $p-\epsilon $, and where $L_{0}v_{i}$ is bounded by $\frac{
2D\lambda _{PF}}{\epsilon (\lambda _{PF}-1)}$. However, there are only a
finite number of words of this length, hence only a finite number of
possible repetition vectors $v_{i}$.\hfill $\square $

\medskip

As defined, the degree of a repetition vector depends on the length vector
of the tiling space. This difference between tiling spaces disappears in the
limit of long repetition vectors:

\begin{lem}
Let $L$ be decomposed as $L= cL_0 + L_r$, where $L_r$ is a linear
combination of (generalized) eigenvectors with eigenvalue less than
$\lambda_{PF}$. For each word $w$, let $r_w$ be the ratio of the
length of $w$ in $\calT_f$ to the length of $w$ in $\calT_0$.  For
each $\epsilon >0$ there exists a length $R$ such that all words with
$|w|>R$ have $|r_w - c| < \epsilon$. \label{lem3.3}
\end{lem}

\medskip 

First consider supertiles of order $n$,
that is words of the
form $\sigma^n a_i$, where $a_i$ is a single letter.  For each such word,
$|r_w-c|$ goes as $(\lambda_2/\lambda_{PF})^n$, where $\lambda_2$ is the 
subleading eigenvalue of $M$. Since the bulk of any long word is made up 
of supertiles of high order, the result follows.
\hfill $\square $
\medskip

\begin{lem}
Let $\epsilon >0$, and let $L$ and $L^{\prime}$ be fixed. There is a
length $ R$ such that, if $v$ is a repetition vector of degree $p$ in
$\calT_f$, representing a word $w$ of length greater than $R$, then
$v$ is a repetition vector of degree $p-\epsilon$ in
$\calT_g$. \label{samerepvecs}
\end{lem}

\medskip \noindent \textbf{\enspace} For the repetition vector $v$
to have degree $p-\epsilon$ in $\calT_g$, the ratio of the lengths
of $w'$ and $w$ in $\calT_g$ must be at least $p-\epsilon$. However,
by Lemma \ref{lem3.3}, 
the ratios of lengths in $\calT_f$ and $\calT_g$ both agree, up to a
small error, with the ratios of lengths in $\calT_0$. 
Since the ratio is at least $p$
in $\calT_f$, by choosing $R$ large enough we can make the ratio at
least $p-\epsilon$ in $\calT_g$. 
\hfill $\square $\medskip

Theorem \ref{families} showed there there are only a finite number of
families of repetition vectors of a given degree.  Lemma
\ref{samerepvecs} shows that the families are essentially the same for
$\calT_f$ and $\calT_g$, up to a small error in the degree.  We next
choose a degree $p$ for which the families are exactly the same for
all choices of tile length, and show that the asymptotic length of the
repetition vectors of degree $p$ yields a conjugacy invariant.

Pick $p_0>1$ such that there exist repetition vectors of degree
strictly greater than $p_0$ in $\calT_0$. By Theorem \ref{families},
there are a finite number of vectors $v_i$ such that every repetition
vector of degree $p_0$ is of the form $M^k v_i$.
Now, for each $k,i$, let $p_{k,i}$ be the maximal degree of $M^k v_i$, 
and let $p_i = \lim_{k \to \infty} p_{k,i}$. Note that for
fixed $i$, the sequence $p_{k,i}$ is nondecreasing in $k$, since the
maximal degree of $\sigma v$ is at least the maximal degree of $v$.
Thus the $p_i$'s are a finite set of real numbers, all greater than
or equal to $p_0$ and some strictly bigger than $p_0$. There are 
real numbers $p$ and $\epsilon$ such that $p > p_0 +
\epsilon$ and such that each $p_i$ is either greater than $p+ \epsilon$
or equal to $p_0$. We restrict our attention to the (nonempty!) 
indices for which $p_i$ is greater than $p+\epsilon$.  Let
$L$, $L^{\prime}$ and $L_0$ be given. By Lemma \ref{samerepvecs},
there is a length $R$ such that the set of repetition vectors $v$ of
degree $p$ and with $L_0 v>R$ is precisely the same for the three
tiling systems $\calT_f$, $\calT_g$ and $\calT_0$. Pick generators
$v_1, \ldots, v_N$ of the finite families of these repetition vectors,
with $L_0 v_i \le L_0 v_{i+1}$, and with $L_0 v_N < \lambda_{PF} L_0
v_1$.

\begin{thm} Let $p>1$ be chosen as in the previous paragraph.
Suppose that $\calT_{f}$ and $\calT_{g}$ are conjugate and that
$\left\{ v_{1},\dots v_{N}\right\} $ is a collection of repetition
vectors of degree $ p$ that generates all repetition vectors of degree
$p.$ Then, given $\delta >0$, for every $i\in \{1,\ldots ,N\}$ and for
any sufficiently large integer $m$, there exist $j,m'$ such
that $|LM^{m}v_{i}-L'M^{m'}v_{j}|<\delta
$. \label{necessary}
\end{thm}

\medskip \noindent \textbf{Proof.\enspace} Let $\phi
:\calT_{f}\rightarrow \calT_{g}$ be the conjugacy. Since $\calT_{f}$
and $\calT_{g}$ are compact metric spaces, $\phi $ is uniformly
continuous. Now suppose $v_{i}$ is a repetition vector of degree
$p$. Let $t_{m,i}=LM^{m}v_{i}$. Then for every tiling $x\in \calT_{f}$
there is a range of real numbers $r$ of size roughly $(p-1)t_{m,i}$,
such that $T_{r}(x)$ is very close to $T_{r+t_{m,i}}(x)$, where
$T_{r}$ denotes translation by $r$. Specifically, as $m\rightarrow
\infty $, the distance from $T_{r+t_{m,i}}(x)$ to $T_{r}(x)$ can be
made arbitrarily small, and the size of the range of $r$'s, divided by
$t_{m,i}$, can be made arbitrarily close to $p-1$. This implies that
$\phi (T_{r}(x))$ and $\phi (T_{r+t_{m,i}}(x))$ are extremely
close. But $\phi (T_{r}(x))=T_{r}(\phi (x))$ and $\phi
(T_{r+t_{m,i}}(x))=T_{r+t_{m,i}}(\phi (x))$, so $\calT_{g}$ admits a
repetition vector whose pairing with $ L'$ is extremely close
to $t_{m,i}$, and whose degree is extremely close to $p$. However, for
$m$ large, all such repetition vectors are of the form $M^{m^{\prime
}}v_{j}$. \hfill $\square $\medskip

\medskip

\noindent {\bf Examples}

\begin{enumerate}

\item {} The
substitution $ a\rightarrow aaaabb$, $b\rightarrow babbba$ has
matrix $ \left(\begin{smallmatrix} 4 & 2\cr2 & 4
\end{smallmatrix}\right) $, whose eigenvalues are 6 and 2.
It is easy to see that $v=(1,0)^{T}$ is a
repetition vector of degree 5, and this vector generates the only
family of repetition vectors with $p=5$. If $\calT _{f}$ and
$\calT_{g}$ are conjugate, for each large $m$ there must exist $
m'$ such that $(LM^{m}-L'M^{m'})v$ is
small.  Since $v$ is full, and since both eigenvalues are large,
this implies there is an integer $k$ such
that $L'=LM^{k}$. 

\item {} The substitution $a\rightarrow aaaabb$,
$ b\rightarrow bbbbaa$, has exactly the same substitution matrix as
the first example,
but the repetitions are different. $v_{1}=(1,0)^{T}$ and
$v_{2}=(0,1)^{T}$ are both repetition vectors with degree 6. The
possibility of having $j\neq i$ in Theorem \ref{necessary} allows for
additional conjugacies. Indeed, $\calT_{f}$
and $\calT_{g}$ are conjugate if (and only if) either $(L'_1, L'_2)=
(L_{1},L_{2})M^{k}$ or $(L'_1,L'_2)=(L_{2},L_{1})M^{k}$ for some
(possibly negative) integer $k$.

\item {} The substitution $a \to b$, $b \to ac$, $c \to ab$, has
a matrix whose eigenvalues are $\tau$, $-1$ and
$-1/\tau$, where $\tau$ is the golden mean.  However, the eigenvalue $-1$
is irrelevant, as its left-eigenvector $(1,-1,0)$ is orthogonal to all
recurrence vectors, and in particular to all repetition
vectors. (Every $a$ is preceded by a $b$, and every $b$ is followed by
an $a$).  This system behaves essentially like the Fibonacci
substitution, with all choices of tile lengths giving conjugate
dynamics, up to an overall scale.

\end{enumerate}

Notice that in the first example, our necessary condition was the same
as the sufficient condition of Theorem \ref{sufficient}, while in the
second and third examples there were conjugacies between length
choices that did not meet the condition of Theorem \ref{sufficient}.
The difference between the first two examples also illustrates that one
cannot obtain sharp conditions on conjugacy from the substitution
matrix alone. One needs some details about the substitution itself.
Our final example further illustrates this point:

\begin{enumerate}

\item[(4)] {} The substitution $a\rightarrow aabaabbba,$
$ b\rightarrow bbabbaaab$ has matrix $M= \left(\begin{smallmatrix} 5 &
4\cr4 & 5
\end{smallmatrix}\right)
=
\left(\begin{smallmatrix}
2 & 1\cr1 & 2
\end{smallmatrix}\right)^{2}$ 
and is the square of the substitution $a\rightarrow aab$, 
$b\rightarrow bba$, so its associated tiling space admits
(among others) conjugacies with $L'=L
\left(\begin{smallmatrix}
2 & 1\cr1 & 2
\end{smallmatrix}\right)^m$, 
where $m$ can be any integer. 
However, the substitution $a\rightarrow aaaaabbbb,$ $b\rightarrow
abbbbbaaa $, with the same matrix $M$, has only one family of repetition
vectors of degree $8$, namely those generated by $\left( 1,0\right) ^{T}$. 
In this case all conjugacies are of the form $L'=LM^{k}
=L \left(\begin{smallmatrix} 2 & 1 \cr 1 & 2 \end{smallmatrix}\right)^{2k}.$

\end{enumerate}

\bigskip
We thank Felipe Voloch, Bob Williams and Charles Radin for
helpful discussions.  This work was partially supported by Texas ARP
grant 003658-152 and by a Faculty Research Grant from the University of 
North Texas.

\end{document}